\documentclass[12pt,thmsa]{article}
\usepackage{amsfonts}
\usepackage{amssymb}
\newtheorem{teo}{Theorem}[section]

\newtheorem{obs2}[teo]{Remark}

\newtheorem{tea}{Theorem}[subsection]

\newtheorem{no2}[teo]{Note}

\newtheorem{no3}[tea]{Note}

\newcommand{\Gal}{{\rm Gal}}

\newcommand{\Ind}{{\rm Ind}}

\newcommand{\Q}{{\mathbb{Q}}}

\begin{document}
\title{{The level $1$ case of Serre's  conjecture revisited
}}
\author{Luis Dieulefait\thanks{Research supported by project MTM2006-04895, MECD, Spain}
\\
Dept. d'Algebra i Geometria, Universitat de Barcelona;\\
Gran Via de les Corts Catalanes 585;
08007 - Barcelona; Spain.\\
e-mail: ldieulefait@ub.edu\\
 }
\date{\empty}

\maketitle

\vskip -20mm

\begin{abstract} We prove existence of conjugate Galois representations, and we use
it to derive a simple method of weight reduction. As a consequence,
an alternative proof of the level 1 case of Serre's conjecture follows.

\end{abstract}

\section{A letter with the results}

$\quad \qquad\qquad\qquad\qquad\qquad\qquad\qquad\qquad\qquad\qquad$Barcelona, April 21, 2007\\

Dear Colleagues:\\

I think there is a simpler way of proving the level $1$ case of Serre's conjecture
and arbitrary weight (i.e., Khare's result). The first steps are of
course as before: you start by proving it for $k = 2$ as in my first work on
Serre's conjecture, and also observing that by Schoof's modularity results
for semistable abelian varieties of conductors $3, 5, 7, 11$ and $13$, you have the
cases of $k = 4, 6, 8, 12, 14$ also covered. These are thus the ``base cases" for
the induction.\\

\noindent I have a procedure to do induction on the weight $k$. So if the representation
has level $1$ and weight (which is thus even) $k > 14$ or $k = 10$, the goal
of the induction step is to reduce such a case of Serre's conjecture to another
case (always with level $1$) of weight $k' < k$.\\

\noindent The setup is an in my first work on Serre's conjecture and the similar work by
Khare-Wintenberger: you use existence of minimal lifts, existence of compatible
families and modularity lifting theorems \`{a} la Wiles to propagate modularity.
We also use existence of weight $2$ lifts, as in Khare's proof (for
simplicity, we remove the distinction between the ordinary and non-ordinary
cases because we now know that weight $2$ lifts exist in both cases). But we
will not use the links that appear in Khare's proof, where he uses an odd
divisor of $p-1$ (for a non-Fermat prime $p$), takes there a non-minimal lift,
thus linking with another compatible family and finally showing that one can
force the weight to decrease.\\

\noindent Thus, what is the new argument? When we have to prove the conjecture
for certain weight $k$, since we can ``switch the prime" we can choose the
characteristic $p$ to be any prime greater or equal to $k-1$. In our proof we
will always choose this suitable prime $p$ to be LARGER that $k$ (i.e., we never
choose $p = k-1$), so the weight two lift will be a lift ``with nebentypus".
For such a $p$-adic Galois representation with nebentypus, let us call it $\rho$, we
want to consider a related one, namely, a conjugate Galois representation $\rho^{\gamma}$.
 The definition of such a representation, (I will prove in the following
lines that it exists!) is as in the case of modular forms, when you change
the Galois representations attached to $f$ by those attached to $f^\gamma$
 where $\gamma$
 is
a field immersion of $K_f$ into the complex numbers, where $K_f$ is the field of
coefficients of $f$ (i.e., an element of the Galois group of the normal closure of
$K_f$ over the rationals).\\
Recall that since $\rho$ is potentially modular (all $p$-adic representations that
appear in our proof are known to be so) it has a ``field of coefficients", namely,
a number field $K$ such that all traces of Frobenius are in this field and they
generate it. We fix an immersion of the algebraic closure of $\Q$ into the one
of $\Q_p$.\\
Using potential modularity, we can show that there exists a ``conjugate" representation
$\rho^\gamma$, where $\gamma$
 is a field immersion of $K$ into the complex, and $\rho^\gamma$
 will
be another potentially modular $p$-adic representation with field of coefficients
$K^\gamma$
 and its traces are $a_p^\gamma$ where $a_p$ are the traces of $\rho$, and its determinant
is also conjugated to the one of $\rho$. The proof is done by imitating the arguments
in the proof of existence of compatible families for potentially modular
Galois representations: you use Brauer's theorem to relate $\rho$ to modular representations
attached to Hilbert modular forms $h_i$ over several real field $F_i$,
then you define in the same way your virtual representation $\rho^\gamma$
 by taking the
same formula except that you replace for each $i$ the $p$-adic representation
attached to $h_i$ by the one attached to $h_i^\gamma$ (and you also change  $\psi$ by  $\psi^\gamma$ 
 for
any character $\psi$  appearing in the formula) and when you want to check that
this is a true Galois representation then the formulas that you have to check,
by just applying the inverse of $\gamma$ 
 to both sides, become equalities that you
know to be true because $\rho$ is a true Galois representation.\\
Thus you conclude the existence of a ``conjugate" Galois representation, and
the local properties of it can be also ``read off" from the modular forms $h_i^\gamma$,
as in the proof of existence of compatible families\footnote{see section 2.1 for more details on this proof}.\\

\noindent I want to consider this conjugate representation only in the situation that it
is useful for the proof of level $1$ Serre's conjecture (after taking a weight two
lift, starting with $p > k$): I have a $p$-adic representation $\rho$ which is potentially
Barsotti-Tate at $p$, unramified outside $p$, and the determinant is $\mu \chi$ where $\chi$ 
is the $p$-adic cyclotomic character and $\mu$ is $\omega^{k-2}$, where $\omega$ is the Teichmuller
lift of the mod $p$ cyclotomic character.\\
Observe that the field of coefficients $K$ of $\rho$ contains the field generated by
the values of $\mu$ (which are roots of unity), let us call $C$ this abelian field
contained in $K$. We take $\gamma$
 such that it acts nontrivially on the roots of
unity that generate $C$, then the representations $\rho$ and $\rho^\gamma$
 will have different
``nebentypus" (i.e., different determinant). For the representation $\rho$ that we
are considering, we know that the inertial Weil-Deligne parameter at $p$ is
exactly:
$$ (\omega^{k-2} \oplus 1, 0)$$
When we proved above the existence of the conjugate representation, we remarked
that the local properties of it can also be deduced from the ones of $\rho$
using potential modularity. In particular, $\rho^\gamma$
 will have nebentypus $\mu^\gamma$
 and it
will be potentially Barsotti-Tate at $p$ with local parameter: $(\mu^\gamma \oplus 1, 0)$ (this
follows from potential modularity, where the field $F$ of modularity can be
taken to be unramified at $p$: over $F$ the representation $\rho$ corresponds to a
Hilbert modular form $h$ and the one we have constructed, $\rho^\gamma$, obviously agrees
with the one corresponding to $h^\gamma$). The local properties at other primes are
proved with the same argument used in the proof of existence of families, in 
our case we conclude that $\rho^\gamma$
 is unramified outside $p$.\\
The idea, as we will see later in more detail, is that by considering this conjugate
representation we can change the nebentypus at $p$: 
 $\gamma$ acts on the roots
of unity in the image of $\mu$ as ``raising to the $i$" for certain exponent $i$, thus
the new inertial parameter at p will be:\\
$(\omega^r \oplus 1, 0)$ for some $r$ (**), we will explain later what values of $r$ are possible
here.\\
We will show that: if $k = 10$ or $k > 14$ we can always take a prime $p > k$
and a suitable $\gamma$
 so that the nebentypus of $\rho^\gamma$
 is as in (**) for an $r$ such
that when we consider the reduction mod $p$ of $\rho^\gamma$, its Serre's weight, which
is known to be (after twisting!) either $r + 2$ or $p + 1 - r$, is, in both cases,
smaller than $k$. Therefore, since it is evident that $\rho$ is modular if and only
if $\rho^\gamma$
 is so, this will conclude the inductive step in the proof of level $1$ Serre
(thanks to modularity lifting theorems).\\
Moreover, in all cases we can just take $p$ to be the smallest prime larger than
$k$, except for $k = 32$ where we need to take $p = 43$.\\

\noindent Remark: If we take $\gamma$
 to be just ``complex conjugation" this is useless. In fact
in this case the Serre's weight of the conjugate (reduced modulo $p$ and after
twisting) will give us again $k$: in fact in the above formula this is one of the
two values we obtain in this case (and using other arguments one can show
that the complex conjugated representation is just a twisted of the given one,
thus the Galois representation has an ``inner twist").\\

\noindent The proof that this procedure always makes $k$ smaller (as long as $k = 10$
or $k > 14$) will be given in two steps: for $k$ up to $36$ we can check it by hand,
something I have already done so I can say which are the values of $p$ (this I
have already said) and which are $\gamma$ 
 and $r$ in each case. In the second step,
for $k > 36$, we work with $p > 37$, and we will use some well-known estimates
on the distribution of primes to prove that the method works, basically we
need to avoid cases like $k = 32$ and $p = 37$ where $p - 1 = 36$ and $k - 2 = 30$
and the ratio here is $36/30 = 6/5$. We want this ratio to be smaller than
that, and it is easy to show that for $p > 37$ it is so.\\

\noindent Let us explain in detail this step: Let $k$ be an even integer with $k = 10$
or $k > 14$, and let $p$ be the smallest prime greater than $k$, except if $k = 32$
where we take $p = 43$.\\
We start with a mod $p$ representation of Serre's weight $k$ and we consider
a weight 2 lift $\rho$. The nebentypus is  $\mu = \omega^{k-2}$ and $\omega$ has order $p - 1$ and
ramifies at $p$ only.\\
Let us call $d = (p-1, k-2)$ and let $m$ be such that $m \cdot d = p-1$. Thus, the
character $\mu$ has order $m$. We choose $\gamma$
 so that the nebentypus is changed to
$\mu^\gamma = \omega^{dt}$ for some $t < m$ with $t$ relatively prime to $m$. We consider $\rho^\gamma$. The
residual Serre's weight of it (after twisting) is equal to $dt+2$ or to $p+1-dt$. Since we want
to CHANGE the Serre's weight after this procedure, we need that the new
nebentypus $\omega^{dt}$ is not equal to $\mu$ nor to the complex conjugate of $\mu$. Thus
we need  $m$ to be such that there are more than $2$ values of $t$, i.e., that for
Euler's $\phi$ function it holds:\\
$\phi(m) > 2$. We will see that we will always have $m > 6$ (or $m = 5$ if $k = 10$,
$p = 11$), so this is true.\\
For the moment, just assume that $m > 6$ (or $m = 5$ if $k = 10$) and we take
the following value for $t$:\\

$\bullet$  $t = (m + 1)/2$ if $m$ is odd\\

$\bullet$ $t = m/2 + 2$ if $m$ is even but not divisible by $4$\\

$\bullet$ $t = m/2 + 1$ if $m$ is divisible by $4$ \\

\noindent Observe that $t$ is always relatively prime to $m$.\\
Let us check, by hand, that for $k$ up to $36$, after taking this conjugate representation,
the residual representation will have a smaller Serre's weight (let
us call $k'$ the Serre's weight after taking Galois conjugation, also recall that
we choose $p = 43$ for $k = 32$):\\

$\bullet$ $k = 10, p = 11$: $d = 2, m = 5, t = 3, dt = 6$; thus: $k' = 8$ or $6$.\\

$\bullet$ $k = 16, p = 17$: $d = 2, m = 8, t = 5, dt = 10$; thus: $k' = 12$ or $8$.\\

$\bullet$ $k = 18, p = 19$: $d = 2, m = 9, t = 5, dt = 10$; thus: $k' = 12$ or $10$.\\

$\bullet$ $k = 20, p = 23$: $d = 2, m = 11, t = 6, dt = 12$; thus $k' = 14$ or $12$.\\

$\bullet$ $k = 22, p = 23$: $d = 2, m = 11, t = 6, dt = 12$; thus $k' = 14$ or $12$.\\

$\bullet$ $k = 24, p = 29$: $d = 2, m = 14, t = 9, dt = 18$; thus $k' = 20$ or $12$.\\

$\bullet$ $k = 26, p = 29$: $d = 4, m = 7, t = 4, dt = 16$; thus $k' = 18$ or $14$.\\

$\bullet$ $k = 28, p = 29$: $d = 2, m = 14, t = 9, dt = 18$; thus $k' = 20$ or $12$.\\

$\bullet$ $k = 30, p = 31$: $d = 2, m = 15, t = 8, dt = 16$; thus $k' = 18$ or $16$.\\

$\bullet$  $k = 32, p = 43$: $d = 6, m = 7, t = 4, dt = 24$; thus $k' = 26$ or $20$.\\

$\bullet$  $k = 34, p = 37$: ...............................................$k' = 22$ or $18$.\\

$\bullet$  $k = 36, p = 37$: ...............................................$k' = 22$ or $16$.\\

\noindent Now we prove the same for $k > 36$ and $p$ the smallest prime larger than $k$
(in particular, $p > 37$). The fact that at the end $k'$ will be smaller than $k$
is based on the fact that two consecutive primes $p_n$ and $p_{n+1}$ are very close (in relative value)
if $p_{n+1} > 37$. We use the same kind of estimates that appear in Khare's
paper, in particular we use the fact that for $x > 100000$ we have Chebyshev's
inequalities for the prime counting function with $A = 1$ and $B = 1.130289$.\\
From this, an elementary argument used also by Khare gives (we need to
take a constant $a > B/A$ and we take $a = 1.144$): For $p_n > 100000$ (the initial
value has not changed because the constant $a$ and $B/A$ are not extremely close\footnote{see section 2.2 for details}), the quotient $p_{n+1}/p_n$ is smaller than
$1.144$.\\
With the help of a computer, we check that in fact this is also true for
$100000 > p_{n+1} > 37$. Thus, if $p_{n+1} > 37$:
$$   p_{n+1}/p_n < 1.144  $$
An obvious corollary of this inequality is the following: For $p_{n+1} > 37$:
$$ (p_{n+1} - 1)/(p_n - 1) < 1.15  $$
In what follows, we will use these two inequalities that hold for $p_{n+1} > 37$.\\
We have a weight $k > 36$ and it is between two primes: $p_n < k < p_{n+1}$, thus
$k - 2  \geq p_n - 1$. The prime $p_{n+1}$ is thus equal to our prime $p$. Then:
   $$ (p_{n+1} - 1)/(k - 2) \leq (p_{n+1} - 1)/(p_n - 1) < 1.15 < 1.2 = 6/5 $$
This implies that $m > 6$. Then we take $t$ as defined before, a value that
tends to half of $m$. An easy computation (see (*) below, where we use
$m > 6$) shows that for such a $t$, if $p_{n+1} > 37$, for the two possible values of
$k'$ that one obtains it always holds: $p_{n+1}/k' > 1.144$.\\
In particular, because of the first inequality for consecutive primes, $k' < p_n$,
therefore $k' < k$ and we are done. This concludes the induction and the new
proof of the level 1 case of Serre's conjecture.\\

\noindent Computation (*): For each of the three cases in the definition of $t$ we take the larger of the two values of $k'$, which is equal to $dt + 2$, and when comparing $p = p_{n+1}$ with $k'$ we obtain the quotients:\\

$\bullet$ $p/k' = (7d+1)/(4d+2)$ or $(9d+1)/(5d+2)$ or $(11d+1)/(6d+2)$......\\

$\bullet$ $p/k' = (10d+1)/(7d+2)$ or $(14d+1)/(9d+2)$ or $(18d+1)/(11d+2)$......\\

$\bullet$ $p/k' = (8d+1)/(5d+2)$ or $(12d+1)/(7d+2)$ or $(16d+1)/(9d+2)$......\\

\noindent In all these cases we easily see that it holds $p/k' > 1.144$. The same holds in the three cases, a fortiori, if we take the smaller of the two possible values of $k'$.\\

\noindent Your comments are suggestions will be strongly appreciated.\\

\noindent Best regards,\\

\noindent Luis Dieulefait

\section{Details}

\subsection{Details on the proof of existence of conjugates}

Let us include, following an editor's suggestion, a more detailed proof of existence of the conjugate representation:\\
As in the proof of existence of compatible families, we start with the relation given by Brauer's formula: Let $F$ be the totally real Galois number field such that, by Taylor's result, we know that the restriction of $\rho$ to it is modular, corresponding to a Hilbert modular form $h$ of parallel weight $2$. We know that $p$ is unramified in $F/\Q$. Let us call $F_i$ the subfields of $F$ such that $\Gal(F/F_i)$ is a solvable group, so by solvable base change we know that over each $F_i$ the restriction of $\rho$ is also modular, corresponding to some Hilbert modular form $h_i$ of parallel weight $2$. Then we have:

$$ \rho = \sum_i  n_i \;  \Ind_{\Gal(F/F_i)}^{\Gal(F/\Q)} \; \rho_{h_i, p} \otimes \phi_i $$
for some characters $\phi_i : \Gal(F/F_i) \rightarrow \bar{\Q}^*$ and integers $n_i$. Observe that here we have used modularity over each $F_i$ to identify the restriction of $\rho$ to each such field with the $p$-adic representation attached to the modular form $h_i$: the key point is that this allows us to consider, for any Galois conjugation $\gamma$ the conjugated representations $\rho_{h_i,p}^\gamma$, equal by definition to the representation $\rho_{h_i^\gamma, p}$ attached to the Hilbert modular form $h_i^\gamma$. Thus, we define as a virtual representation:

$$ \rho^\gamma = \sum_i  n_i  \; \Ind_{\Gal(F/F_i)}^{\Gal(F/\Q)} \; \rho_{h_i, p}^\gamma \otimes \phi_i^\gamma $$

To check that it is a true Galois representation, we proceed as in the proof of existence of compatible families given in [Di1] and we compute the inner product
 $(\rho^\gamma, \rho^\gamma)$, via an application of Frobenius reciprocity and Mackey's formula (cf. [Ta3], section 5.3.3) we obtain:
 
 $$(\rho^\gamma, \rho^\gamma) = \sum_{i,j}  \; \sum_{g \in G_{F_i} \backslash  G_\Q  / G_{F_j} }    t_{i,j,g}$$
 where $G_K$ denotes the absolute Galois group of $K$ for any number field $K$, and $t_{i,j,g}$ is defined as follows:
$t_{i,j,g} = n_i \cdot n_j$ if
$$  \rho_{h_i, p}^\gamma \otimes \phi_i^\gamma |_{G_{F_i \; {^gF}_j}}  \cong 
 c_g \circ \rho_{h_j, p}^\gamma \otimes \phi_j^\gamma|_{G_{^{g^{-1}}F_i \; F_j         }       } $$
 and $t_{i,j,g} = 0$ otherwise. In the above formula, $c_g$  transforms, for $K = {^{g^{-1}}F_i \; F_j}$, representations of $G_K$ into representations of
  $G_{^gK}$ by conjugation, i.e., transforms $\sigma$ into $\sigma(g^{-1}-g)$.\\
 Thus it is easy to see that the value of this inner product is the same as that of the inner product $(\rho, \rho)$ just by the following elementary and fundamental principle: 
 $$ A = B \Leftrightarrow A^\gamma = B^\gamma $$ 
 for any pair of algebraic numbers $A, B$ and any Galois conjugation $\gamma$. \\
 Thus $(\rho^\gamma, \rho^\gamma) = (\rho, \rho) = 1$, the last equality follows from the fact that $\rho$ is a true, irreducible, Galois representation. Then we conclude that $\rho^\gamma$ is also a true, irreducible, Galois representation, and this concludes the proof since by construction it is clear that $\rho^\gamma$ satisfies the definition of ``conjugate" representation that we have given in the previous section. \\

\subsection{On the quotient of consecutive primes}

Starting from the following Chebyshev's inequalities for the prime counting function:
 
 $$A \;  \frac{x}{\log x} < \pi(x) < B \;  \frac{x}{\log x} $$
 with $A=1$ and $B= 1.130289$ which are known to hold for any $x > x_0 = 100000$, as in Khare's paper if we take $a > C := B/A = B$ then we also have:
 $p_{n+1} / p_n < a$ for any\footnote{in [Kh], due to a small typo, the value $a^{C/(a-C)}$ appears as $aC/(a-C)$}
 $p_n > \max(x_0, a^{C/(a-C)})$. We have chosen $a= 1.144$ and since for this value we easily check that  $a^{C/(a-C)} =65530.89... <100000$ then we conclude that for $p_n > 100000$ it holds $p_{n+1}/p_n < 1.144$. \\

\section{Bibliography}

All the technical details needed to make the above proofs work have already
been proved in previous papers on Serre's conjecture by the author and by
Khare-Wintenberger or Khare. In particular, in these papers the following
results are proved: existence of compatible families, existence of minimal
lifts and lifts with prescribed local properties, and an explanation that available
modularity lifting results can be applied to the $p$-adic representations
appearing in the above proof (both in the residually reducible and in the
residually modular case), thus allowing to propagate modularity whenever
needed. Also the base cases of small weight are proved in these previous
works. 
Thus we just indicate in this bibliography these papers, together
with Taylor's papers on potential modularity and applications which constitute
the main tool used in all these works (and used in particular in the proof
of existence of conjugate Galois representations given in this note), but the
reader is advised to consult also the references therein:\\

\noindent [Di1] Dieulefait, L., {\it Existence of compatible families and new cases of the Fontaine-Mazur conjecture},
 J. Reine Angew. Math. {\bf 577} (2004) 147-152 
 \newline
 [Di2] Dieulefait, L., {\it The level $1$ weight $2$ case of Serre's conjecture}, preprint (2004); published: December 2007
\newline
[Di3] Dieulefait, L., {\it Remarks on Serre's modularity conjecture}, preprint (2006)
\newline
[Kh] Khare, C., {\it Serre's modularity conjecture: the level $1$ case}, Duke Math. J. {\bf 134} (2006)  557-589
\newline
[K-W1] Khare, C., Wintenberger, J.-P., {\it On Serre's conjecture for $2$-dimensional mod $p$ representations of the Galois group of $\Q$}, preprint (2004)
\newline
[K-W2] Khare, C., Wintenberger, J.-P., {\it Serre's modularity conjecture (1)}, preprint (2006)
\newline
[Ta1] Taylor, R., {\it Remarks on a conjecture of Fontaine and Mazur},
J. Inst. Math. Jussieu {\bf 1} (2002)
\newline
[Ta2] Taylor, R., {\it On the meromorphic continuation of degree two
 L-functions}, Documenta Mathematica, Extra Volume: John Coates' Sixtieth Birthday (2006) 729-779 
\newline
[Ta3] Taylor, R., {\it Galois representations}, Annales Fac. Sc. Toulouse {\bf 13} (2004) 73-119

\end{document}